\documentclass[preprint,3p,12pt]{elsarticle}

\usepackage{amssymb}
\usepackage{latexsym}

\usepackage{url}
\usepackage{xcolor}
\definecolor{newcolor}{rgb}{.8,.349,.1}

\usepackage{bm}
\usepackage[ruled,vlined]{algorithm2e}
\usepackage[caption = false]{subfig}
\usepackage{amsmath,amsfonts}
\usepackage{hyperref}
\hypersetup{colorlinks=true}
\usepackage{capt-of}
\usepackage{array}
\usepackage{stackengine}
\usepackage[normalem]{ulem}

\newcommand{\dv}[2]{\frac{\partial #1}{\partial #2}}

\graphicspath{{./figs/}}

\begin{document}


\begin{frontmatter}

\title{
	A Variable Eddington Factor Model for Thermal Radiative Transfer with Closure based on Data-Driven Shape Function
}

\author[1,3]{Joseph M. Coale}

\author[2,4]{Dmitriy Y. Anistratov}

\address[1]{Computational Physics and Methods Group, Los Alamos National Laboratory, Los Alamos, NM}
\address[2]{Department of Nuclear Engineering, North Carolina State University, Raleigh, NC}
\address[3]{jmcoale@lanl.gov}
\address[4]{anistratov@ncsu.edu}


\begin{abstract}
A new variable Eddington factor (VEF) model is presented for nonlinear problems of thermal radiative transfer (TRT).
The VEF model is a data-driven one that acts on known (a-priori)  radiation-diffusion solutions for material temperatures in the TRT problem.
A linear  auxiliary  problem is constructed for the radiative transfer equation (RTE) with opacities and emission source evaluated at the known material temperatures.
The solution to this RTE approximates the specific intensity distribution for the problem in all phase-space and time.
It is  applied as a shape function to define the Eddington tensor for the presented VEF model.
The shape function computed via the auxiliary RTE problem will capture some degree of
transport effects within the TRT problem. The VEF moment equations closed with this 
approximate Eddington tensor will thus carry with them these captured transport effects.
In this study, the temperature data comes from multigroup $P_1$, $P_{1/3}$, and  flux-limited  diffusion radiative transfer (RT)  models.
The proposed VEF model can be interpreted as   a transport-corrected diffusion reduced-order model.
Numerical results are presented on the Fleck-Cummings test problem which models a supersonic wavefront of radiation.
The presented VEF model is shown to reliably improve accuracy by 1-2 orders of magnitude compared to the considered radiation-diffusion model solutions to the TRT problem.
\end{abstract}


\end{frontmatter}


%
%
\section{Introduction}
The modeling and simulation of thermal radiation transport is an important consideration for many applications. Such phenomena can be found in a wide array of different fields, including: high-energy density (HED) physics, astrophysics, plasma physics, fire and combustion physics, and atmospheric and ocean sciences \cite{zel-1966,mihalas-FRH-1984,shu-astro,thomas-stamnes-atm,faghri-sunden-2008,drake-hedp}.
Multiphysical models of these phenomena comprise complex systems of partial differential equations (PDEs) which must be solved by means of numerical simulation.
Some challenges are associated with the numerical simulation these systems.
The involved equations are characterized by tight coupling, strong nonlinearities and multiscale behavior in space-time.

Since radiative transfer (RT) is an important mechanism for energy redistribution in these phenomena,
a photon transport model must be included in the systems of equations to model RT effects.
The equations of RT are high-dimensional, their solution typically depending on 7 independent variables in 3D geometry.
Thus along with the aforementioned challenges, a massive number of degrees of freedom (DoF) must be used in calculations to adequately describe the overall solution.
For a simulation that discretizes each independent variable with an $x-$point grid, $x^7$ DoF are required, making it reasonable to reach trillions or quadrillions of DoF.
The resulting computational load and memory occupation can become intractable for large simulations without the use of exascale computing resources.

A common approach to develop models for simulation of RT phenomena in multiphysics problems is to apply methods of dimensionality reduction for the RT component. This will significantly reduce the required computational resources in exchange for some level of approximation and modeling error in the resulting solution.
The  goal for an RT model is to find balance between computational load with the desired fidelity of solutions.
There exist many well-known RT models based on moment equations with approximate closures
which have seen extensive use in application, such as $M_N$ methods that use maximum entropy closure relations, the $P_1$ and $P_{1/3}$ models, and flux-limited diffusion models \cite{olson-auer-hall-2000,morel-2000,simmons-mihalas-2000,hauck-2011,hauck-2012,m1-2017}.
  Variable Eddington factor (VEF) models makeup another such class of RT models \cite{eddington-1926,pomraning-1969,l-p-1981}.
VEF models are constructed by reformulating the radiation pressure tensor in terms of the Eddington tensor which brings closure to the moment equations.
There exist many ways to construct approximation for the Eddington tensor.
Some commonly used VEF models apply Wilson, Kershaw and Levermore closures \cite{wilson-1970,kershaw-1976,levermore-1984,levermore-1996,m1-2017}.
Numerical methods for solving the Boltzmann transport equation have been developed based on the first two angular moment equations with exact closure defined by the Eddington tensor
\cite{gol'din-1964,auer-mihalas-1970,gol'din-1972,PASE-1986,winkler-norman-mihalas-85}.

The anisotropic diffusion (AD) model has been developed for thermal radiative transfer (TRT) and other particle transport applications \cite{trahan-larsen-2009,johnson-larsen-2011,trahan-larsen-2011,johnson-larsen-2011-2}.
A diffusion equation is constructed which is closed by means of an AD coefficient.
The AD coefficient is defined via the solution to an auxiliary transport problem which takes the form of an angular and spatially dependent shape function.
The special shape function accounts for some degree of transport effects in the TRT problem.
This yields the tensor-diffusion moment equations with transport-corrected AD coefficients.
The AD model  uses just one transport sweep to solve for the auxiliary function.

A new approach based on data-driven reduced-order models (ROMs)
has been gaining popularity in recent years which make use of
data-based methodologies
to dimensionality reduction.
Data-driven models have been developed for
(i) linear particle transport problems
\cite{buchan-2015,tencer-2017,soucasse-2019,behne-ragusa-morel-2019,hardy-morel-ahrens-2019,prince-ragusa-2019-1,prince-ragusa-2019-2,Dominesey-2019,peng-mcclarren-frank-2019,peng-mcclarren-tans2019,pozulp-2019,peng-mcclarren-frank-2020,pozulp-brantley-palmer-vujic-2021,peng-mcclarren-2021,behne-ragusa-tano-2021,choi-2021,elhareef-wu-ma-2021,peng-chen-cheng-li-2022,hughes-buchan-2022,huang-I-2022}
(ii) nonlinear RT problems \cite{pinnau-schulze-2007,qian-wang-song-pant-2015,fagiano-2016,alberti-palmer-2018,jc-dya-m&c2019,jc-dya-tans2019,jc-dya-m&c2021,girault-2021,alberti-jqsrt-2022,jc-dya-jqsrt2023,dya-jc-m&c2023-2,jmc-dya-arxiv2023}, and
(iii) various problems in nuclear reactor-physics  \cite{cherezov-sanchez-joo-2018,alberti-palmer-2019,german-ragusa-2019-1,german-ragusa-2019-2,alberti-palmer-2020,german-tano-ragusa-2021,elzohery-roberts-2021-1,elzohery-roberts-2021-2,elzohery-roberts-2021-3,phillips-2021,Dominesey-2022}.
The fundamental idea behind these ROMs is to leverage databases of solutions to their problems of interest (known a-priori) to develop some reduction in the dimensionality for their involved equations.
By nature, these models are problem-dependent since they are formulated using chosen datasets, and this allows for higher levels of accuracy than displayed by other types of ROMs (within the regime of parameters covered by the datasets).

In this paper we present a data-driven VEF   model  for the fundamental  TRT  problem.
This kind of TRT problem models a supersonic flow of radiation through matter
and neglectshydrodynamic motion of the underlying material and heat conduction \cite{Moore-2015}.
Note that the TRT problem is characterized by the same fundamental challenges as the more general class of problems (e.g. radiation-hydrodynamics problems) and serves as a useful computational platform for the development of new models.
There exist several pathways to obtain an approximate Eddington tensor by data-driven means if data on the RT solution
for a subset of problem parameters
can be obtained \cite{jc-dya-m&c2019,jc-dya-m&c2021,jc-dya-jqsrt2023,dya-jc-m&c2023-2,jmc-dya-arxiv2023}.
For the model developed here, the Eddington tensor is computed
from approximate solution data to the TRT problem generated
with radiation-diffusion models.
It has been previously shown that the   solution of the Boltzmann transport equation  computed with a scattering source term evaluated by a diffusion solution yields a sufficiently accurate shape function for estimation of the Eddington tensor in linear problems \cite{dya-ns-2012,ns-dya-mla-2014}.
An extension of this idea to the nonlinear TRT problem is to use the material temperatures evaluated with a radiation-diffusion model to compute the opacity and emission source in the radiative transfer equation (RTE).
This step of  the  new model  uses only one transport sweep to calculate a shape function
accounting approximately for transport effects, and to generate the Eddington tensor for the TRT problem.

The remainder of the paper is as follows. The TRT problem is described in Sec. \ref{sec:trt}, along
with   definitions for  several classical moment-based   RT models 
applied to generate data for computation of the auxiliary shape function. The new data-driven VEF model is formulated in Sec. \ref{sec:vef}. Numerical results are given in Sec. \ref{sec:results}, followed by conclusions in Sec. \ref{sec:conclusion}.

%
%
\section{Thermal Radiative Transfer and Models Based on Moment Equations} \label{sec:trt}
The TRT problem is formulated by the multigroup RTE
\begin{subequations}\label{bte}
	\begin{gather}
		\frac{1}{c}\dv{I_g}{t} + \boldsymbol{\Omega}\cdot\boldsymbol{\nabla}I_g + \varkappa_g(T)I_g = \varkappa_g(T)B_g(T), \\
		I_g|_{\boldsymbol{r}\in\partial\Gamma}=I_g^\text{in}\ \ \text{for}\ \ \boldsymbol{\Omega}\cdot\boldsymbol{n}_\Gamma<0,\quad I_g|_{t=0}=I_g^0, \\
		\boldsymbol{r}\in\Gamma,\quad  t> 0,\quad \boldsymbol{\Omega}\in\mathcal{S},\quad g=1,\dots,G, \nonumber
	\end{gather}
\end{subequations}
and the material energy balance (MEB) equation
\begin{equation} \label{meb}
	\dv{\varepsilon(T)}{t} = {\sum_{g=1}^{G}\bigg(\int_{4\pi}I_gd\Omega - 4\pi B_g(T)\bigg)\varkappa_g(T)},\quad T|_{t=0}=T^0 \, ,
\end{equation}
where
$\boldsymbol{r}$ is spatial position, $t$ is time, $g$ is the frequency group index, $\Gamma$ is the spatial domain, $\partial\Gamma$ is the boundary surface of $\Gamma$,
$\boldsymbol{n}_\Gamma$ is the unit outward normal to $\partial\Gamma$,
$I_g(\boldsymbol{r},\boldsymbol{\Omega},t)$ is the group specific photon intensity,
$T(\boldsymbol{r},t)$ is the material temperature,
$\varkappa_{g}(\boldsymbol{r},t;T)$ is the group material opacity, $\varepsilon(\boldsymbol{r},t;T)$ is the material energy density,  and $B_g(\boldsymbol{r},t;T)$ is the group Planckian function given by
\begin{equation} \label{eq:planck}
	B_g(T) =  \frac{2}{c^2h^2}\int_{\nu_{g-1}}^{\nu_{g}}  \frac{\nu^3 }{e^{\frac{\nu}{T}}-1} d\nu .
\end{equation}
$c$ is the speed of light, $h$ is Planck's constant, $\nu$ is photon frequency.

There are several  TRT models  which apply moment equations
for the group radiation energy density
\begin{equation} \label{Eg}
	E_g =  \frac{1}{c} \int_{4\pi}  {I}_g \ d\Omega ,
\end{equation}
and flux
\begin{equation}  \label{Fg}
	\boldsymbol{F}_g =  \int_{4\pi} \boldsymbol{\Omega} {I}_g \ d\Omega,
\end{equation}
to approximate the RTE.
The $P_1$ model  is defined by the  multigroup $P_1$ equations for radiative transfer, given by
\begin{subequations} \label{eq:p1_eqs}
	\begin{gather}
		\dv{E_g}{t} + \boldsymbol{\nabla}\cdot\boldsymbol{F}_g + c\varkappa_g(T)E_g = 4\pi\varkappa_g(T)B_g(T), \label{p1_ebal}\\
		\frac{1}{c}\dv{\boldsymbol{F}_g}{t} + \frac{c}{3}\boldsymbol{\nabla} E_g + \varkappa_g(T) \boldsymbol{F}_g = 0, \label{p1_m1}\\
		\boldsymbol{F}_g|_{\boldsymbol{r}\in\partial\Gamma} = \frac{1}{2}E_g|_{\boldsymbol{r}\in\partial\Gamma}+2{F}_g^\text{in},\quad  \\
		E_g|_{t=0}=E_g^0,\quad  \boldsymbol{F}_g|_{t=0}=\boldsymbol{F}_g^0 .
	\end{gather}
\end{subequations}
The hyperbolic time-dependent $P_1$ equations are derived from the RTE by taking its
$0^\text{th}$ and $1^\text{st}$ angular moments. Closure for the moment equations is
formulated by defining the highest ($2^\text{nd}$) angular moment
\begin{equation}
	H_g =  \int_{4\pi} \boldsymbol{\Omega}\otimes\boldsymbol{\Omega} {I}_g\ d\Omega
\end{equation}
with the expansion
\begin{equation}
	I_g = \frac{1}{4 \pi} (cE_g + 3 \boldsymbol{\Omega}\cdot \boldsymbol{F}_g) \, .
\end{equation}
This approximation yields
\begin{equation}
	H_g =  \frac{c}{3} E \, .
\end{equation}

The $P_{1/3}$ model for radiative transfer is  formulated
by the balance equation \eqref{p1_ebal} and  the modified first moment equation given by  \cite{olson-auer-hall-2000,morel-2000}
\begin{equation} \label{p1/3-m1}
	\frac{1}{3c}\dv{\boldsymbol{F}_g}{t} + \frac{c}{3}\boldsymbol{\nabla} E_g + \varkappa_g(T) \boldsymbol{F}_g = 0 \, ,
\end{equation}
The  factor   $\frac{1}{3}$  at the time-derivative term in Eq.  \eqref{p1/3-m1}
 produces  the correct the propagation speed of radiation in vacuum.

The flux-limited diffusion (FLD) models are defined by
the time-dependent multigroup  diffusion equations \cite{l-p-1981,krumholz-2007,kuiper-2010,a&a-2015,tetsu-2016}
\begin{subequations}
	\label{eq:diff}
	\begin{gather}
		\label{eq:diff_eq}
		\dv{E_g}{t} + c\boldsymbol{\nabla}\cdot(-D_g\boldsymbol{\nabla} E_g) + c\varkappa_{g}(T)E_g = 4\pi\varkappa_{g}(T)B_g(T), \\
		\boldsymbol{n}_\Gamma\cdot(-cD_g\boldsymbol{\nabla} E_g)|_{\boldsymbol{r}\in\partial\Gamma}=\frac{1}{2}E_g|_{\boldsymbol{r}\in\partial\Gamma}+2{F}_g^\text{in},\quad E_g|_{t=0}=E_g^0,
	\end{gather}
\end{subequations}
and
\begin{equation}
	\boldsymbol{F}_g =  -cD_g\boldsymbol{\nabla} E_g \, ,
\end{equation}
where $D_g$ is the group diffusion coefficient.
In this model, the time derivative of the flux in the first moment equation is neglected.
This leads to a parabolic time-dependent equation for $E_g$ with the diffusion coefficient defined by
\begin{equation} \label{diff-coef}
	D_g  =   \frac{1}{3\varkappa_{g}(T)} \, .
\end{equation}
In general, the solution of the diffusion equation \eqref{eq:diff_eq} with  $D_g$ defined by Eq. \eqref{diff-coef}  does not satisfy  the flux-limiting condition
\begin{equation} \label{f-limit}
	\frac{| \mathbf{F}_g|}{cE_g}   \le 1 \, ,
\end{equation}
stemming from definitions of the radiation density and flux (Eqs.  \eqref{Eg} and \eqref{Fg}).
The FLD models introduce  modifications of the diffusion coefficient  to meet the condition in Eq. \eqref{f-limit}.
In this study,  we consider the coefficient proposed by E. Larsen \cite{olson-auer-hall-2000}
\begin{equation}
	D_g(T,E_g) = \bigg[ \ \big( 3\varkappa_{g}(T) \big)^2 + \bigg( \frac{1}{E_g}\boldsymbol{\nabla}E_g \bigg)^2 \ \bigg]^{\frac{1}{2}}.
\end{equation}

%
%
\section{Variable Eddington Factor Model for TRT with Diffusion-Based Shape Function} \label{sec:vef}
The Variable Eddington  factor  method is defined by the balance equation (Eq. \eqref{p1_ebal}) and the first moment equation
\begin{equation} \label{lovef1}
	\frac{1}{c}\dv{\boldsymbol{F}_g}{t} +c\boldsymbol{\nabla} (\mathfrak{f}_gE_g) + \varkappa_g(T) \boldsymbol{F}_g = 0 \, ,
\end{equation}
where closure is defined with
\begin{equation}
	H_g =  c\mathfrak{f}_g[\tilde I] E_g,
\end{equation}
by means of the Eddington tensor given as
\begin{equation} \label{approx-ET}
	\mathfrak{f}_g[\tilde I]  =  \frac{\int_{4\pi} \boldsymbol{\Omega}\otimes\boldsymbol{\Omega} {\tilde I}_g d\Omega}{\int_{4\pi} {\tilde I}_g d\Omega} \, .
\end{equation}
Here $\tilde I_g$ is an approximation of the photon intensity.
There exist a group of VEF models
which use an approximation of the Eddington tensor defined via the first two moments of the photon intensity \cite{kershaw-1976,minerbo-1978,levermore-1984,su-mgf-ewl-2001,swesty-2009,skinner-2013,klassen-2014,tetsu-2016,Zhang-2013,m1-2017}.

To define the Eddington tensor, we
formulate a model in which the material temperature distribution $\tilde{T}$ for a TRT problem is calculated with one of the radiation-diffusion models described in Sec. \ref{sec:trt}.  A  linear RTE is then
defined by available $\tilde{T}(\boldsymbol{r},t)$ for $t\in[0,t^\text{end}]$ and $\boldsymbol{r}\in\Gamma$
\begin{subequations}
	\label{eq:tcd_bte}
	\begin{gather}
		\frac{1}{c}\dv{\tilde{I}_g}{t} + \boldsymbol{\Omega} \cdot \boldsymbol{\nabla} \tilde{I}_g + \varkappa_{g}(\tilde{T})\tilde{I}_g =  \varkappa_{g}(\tilde{T}) B_g(\tilde{T}),  \label{eq:bte_fixsrc} \\
		\tilde{I}_g |_{\boldsymbol{r}\in\partial\Gamma} = I_g^\text{in}\ \ \text{for}\ \ \boldsymbol{n}_\Gamma\cdot\boldsymbol{\Omega}<0,\quad \tilde{I}_g |_{t=t_0} = I_g^0,\\
		\boldsymbol{r}\in\Gamma,\quad t\in[0,t^\text{end}],\quad \boldsymbol{\Omega}\in\mathcal{S},\quad g=1,\dots,G. \nonumber
	\end{gather}
\end{subequations}
The solution of the auxiliary RTE problem \eqref{eq:tcd_bte} gives an approximate distribution of radiation intensities $\tilde{I}_g$ which
accounts for the transport effects of the TRT problem and can used as a shape function to compute approximate Eddington tensor \eqref{approx-ET}.
The boundary conditions for the VEF moment equations are defined in terms of $\tilde{I}_g$  as follows \cite{gol'din-1964,gol'din-1972}:
\begin{equation}
	\boldsymbol{n}_\Gamma\cdot\boldsymbol{F}_g|_{\boldsymbol{r}\in\partial\Gamma}=cC_g[\tilde{I}_g](E_g|_{\boldsymbol{r}\in\partial\Gamma} - E^\text{in})+{F}_g^\text{in},
	\label{eq:newbc}
\end{equation}
where
\begin{equation}
	C_g[\tilde{I}_g] = \frac{\int_{\boldsymbol{n}_\Gamma\cdot\boldsymbol{\Omega}>0}\boldsymbol{\Omega}\tilde{I}_g\ d\Omega}{\int_{\boldsymbol{n}_\Gamma\cdot\boldsymbol{\Omega}>0}\tilde{I}_gd\Omega}.
\end{equation}
The RTE \eqref{eq:tcd_bte}  with the given function of  temperature can be efficiently solved with a single transport sweep per time step.
To solve the Eqs. \eqref{eq:tcd_bte}
ray tracing techniques (aka the method of long characteristics) are applied \cite{goldin-1960,takeuchi-1969,askew-1972,brough-chudley-1980,askew-roth-1982,casmo-1993,zika-adams-2000,dya-jc-m&c2023}.
In sum,  the data-driven VEF model for TRT is constructed with:
\begin{itemize}
	\item Radiation-diffusion solution data for material temperatures $\tilde{T}$,
	\item The RTE with opacity and Planckian source evaluated with $\tilde{T}$ (Eqs. \eqref{eq:tcd_bte}),
	\item The VEF equations (Eqs. \eqref{p1_ebal} \& \eqref{lovef1}), where the Eddington tensor is defined via Eq. \eqref{approx-ET} and boundary conditions given by Eq. \eqref{eq:newbc}.
\end{itemize}
Hereafter we refer to this model as the data-driven VEF model (DD-VEF).

\begin{algorithm}[ht!]
	\SetAlgoLined
	\vspace*{.1cm}
	Input: $\{\boldsymbol{\tilde{T}}(t^n)\}_{n=1}^N$\\
	$n=0$\\
	\While{$t^n<t^\text{end}$}{
		$n=n+1$\\
		$\boldsymbol{\tilde{T}}(t^n)\leadsto\boldsymbol{\tilde{I}}_g(t^n)$ (Eqs. \eqref{eq:tcd_bte})\\
		$\boldsymbol{\tilde{I}}_g(t^n)\leadsto{\boldsymbol{\tilde{f}}}_g(t^n)$ (Eq. \eqref{approx-ET})
	}
	Output: $\{\ {\boldsymbol{\tilde{f}}}_g(t^n)\}_{n=1}^N$
	\caption{Offline phase of DD-VEF model
		\label{alg:off}}
\end{algorithm}

\begin{algorithm}[ht!]
	\SetAlgoLined
	\vspace*{.1cm}
	Input: $\{{\boldsymbol{\tilde{f}}}_g(t^n)\}_{n=1}^N$\\
	$n=0$\\
	\While{$t^n<t^\text{end}$}{
		$n=n+1$\\
		${\boldsymbol{\tilde{f}}}_g(t^n)\leadsto\{\boldsymbol{T}(t^n),\ \boldsymbol{{E}}_g(t^n),\ \boldsymbol{{\mathcal{F}}}_g(t^n)\}$ (Eqs.  \eqref{p1_ebal}, \eqref{lovef1}, \eqref{eq:newbc})
	}
	Output: $\{\ \boldsymbol{{T}}(t^n),\ \boldsymbol{{E}}_g(t^n),\ \boldsymbol{{\mathcal{F}}}_g(t^n)\}_{n=1}^N$
	\caption{Online phase of DD-VEF model
		\label{alg:on}}
\end{algorithm}

The process of solving TRT problems with the DD-VEF model can be explained as a two-phase methodology, which is outlined in Algorithms \ref{alg:off} and \ref{alg:on}.
The first (offline) phase demonstrated by Algorithm \ref{alg:off} represents the `data-processing' operations to prepare the Eddington tensor closure data. 
The required input is an already known approximate material temperature distribution $\tilde{T}$ for the entire spatial and temporal interval of interest. 
If we define a simulation with $X$ spatial grid cells and $N$ time steps, then this input data is the set $\{\boldsymbol{\tilde{T}}(t^n)\}_{n=1}^N$ where $\boldsymbol{\tilde{T}}(t^n)\in\mathbb{R}^X$.
At each $n^\text{th}$ time step, Eq. \eqref{eq:tcd_bte} is solved using $\boldsymbol{\tilde{T}}(t^n)$ for the vector of discrete radiation intensities in phase space $\boldsymbol{\tilde{I}}_g$, which give rise to the approximate Eddington tensor on the discrete grid $\boldsymbol{\tilde{f}}_g(t^n)$ via Eq. \eqref{approx-ET}.
The discrete Eddington tensor values at each time step can be collected and stored in a dataset for later use with the DD-VEF model.
This process of preparing the Eddington tensor data is referred to as the offline-phase because it only must completed once per temperature distribution $\tilde{T}$, and the calculated Eddington tensor data can be stored away for later use.

The second (online) phase is outlined in Algorithm \ref{alg:on} and represents the operations required to solve a given TRT problem with the DD-VEF model. Taking as input the Eddington tensor data calculated in Algorithm \ref{alg:off}, The DD-VEF equations are solved at each time step to generate vectors for temperature $\boldsymbol{{T}}(t^n)$, radiation energy densities $\boldsymbol{{E}}_g(t^n)$ and radiation fluxes $\boldsymbol{{\mathcal{F}}}_g(t^n)$ over all phase space. In this configuration of offline/online phases, only Algorithm \ref{alg:on} must be completed for any given TRT simulation, assuming Algorithm \ref{alg:off} was completed some time in the past to prepare the required datasets. It is important to note however, that both phases can be combined to save on storage requirements. In this case given the input for $\{\boldsymbol{\tilde{T}}(t^n)\}_{n=1}^N$, at each time step the approximate Eddington tensor is calculated and immediately used with the DD-VEF equations to generate the TRT solution for the time step.

%
%
\section{Numerical Results} \label{sec:results}
The DD-VEF model is analyzed with numerical testing on the classical Fleck-Cummings (F-C) test \cite{fleck-1971} in 2D Cartesian geometry.
This test takes the form of a homogeneous square domain with sides 6 cm in length, whose material is defined with spectral opacity
\begin{equation}
	\varkappa_\nu = \frac{27}{\nu^3}(1-e^{-\nu/T}).
\end{equation}

\noindent Here $\nu$ and $T$ are measured in KeV.
The left boundary is subject to an isotropic, black-body radiation source at a temperature of $T^\text{in}=1$ KeV. All other boundaries are vacuum.
The initial temperature of the domain is $T^0=1$ eV. The material energy density of the material is a linear one $\varepsilon=c_vT$ where $c_v=0.5917 a_R (T^\text{in})^3$.
The problem is solved on the interval $t\in[0,6\, \text{ns}]$ with 300 uniform time steps $\Delta t = 2\times 10^{-2}$ ns. The phase space is discretized using a $20\times 20$ uniform orthogonal spatial grid, 17 frequency groups (see Table \ref{tab:freq_grps}) and 144 discrete directions. The Abu-Shumays angular quadrature set is used \cite{abu-shumays-2001}. The implicit backward-Euler time integration scheme is used to discretize all equations in time. The BTE is discretized in space with the method of conservative long characteristics \cite{dya-jc-m&c2023}, and all low-order equations use a second-order finite-volumes scheme \cite{pg-dya-jcp-2020}.

Note the full-order model (FOM) for this TRT problem is formulated as the RTE coupled with the MEB.
Three radiation diffusion models are considered to generate $\tilde{T}$: multigroup FLD, $P_1$ and $P_{1/3}$ (see Sec. \ref{sec:trt}). The physics embedded in $\tilde{T}$ will vary depending on which diffusion type model is used in its computation. For instance the FLD, $P_1$ and $P_{1/3}$ models may all produce different propagation speeds (and spectral distributions) of the radiation wavefront \cite{olson-auer-hall-2000}. These effects change how energy is redistributed in the F-C test and alters the distribution of material temperatures in space-time.

\begin{table*}[ht!]
	\centering
	\caption{Upper boundaries for each frequency group}
	\label{tab:freq_grps}
	\begin{tabular}{|l|l|l|l|l|l|l|l|l|l|}
		\hline
		$g$ &1&2&3&4&5&6&7&8&9 \\ \hline
		$\nu_{g}$ [KeV]
		& 0.7075
		& 1.415
		& 2.123
		& 2.830
		& 3.538
		& 4.245
		& 5.129
		& 6.014
		& 6.898 \\ \hline\hline
		$g$&10&11&12&13&14&15&16&17& \\ \hline
		$\nu_{g}$ [KeV]
		& 7.783
		& 8.667
		& 9.551
		& 10.44
		& 11.32
		& 12.20
		& 13.09
		& 1$\times 10^{7}$ & \\ \hline
	\end{tabular}
\end{table*}

\begin{figure}[ht!]
	\centering
	\subfloat[Material Temperature]{\includegraphics[width=.5\textwidth]{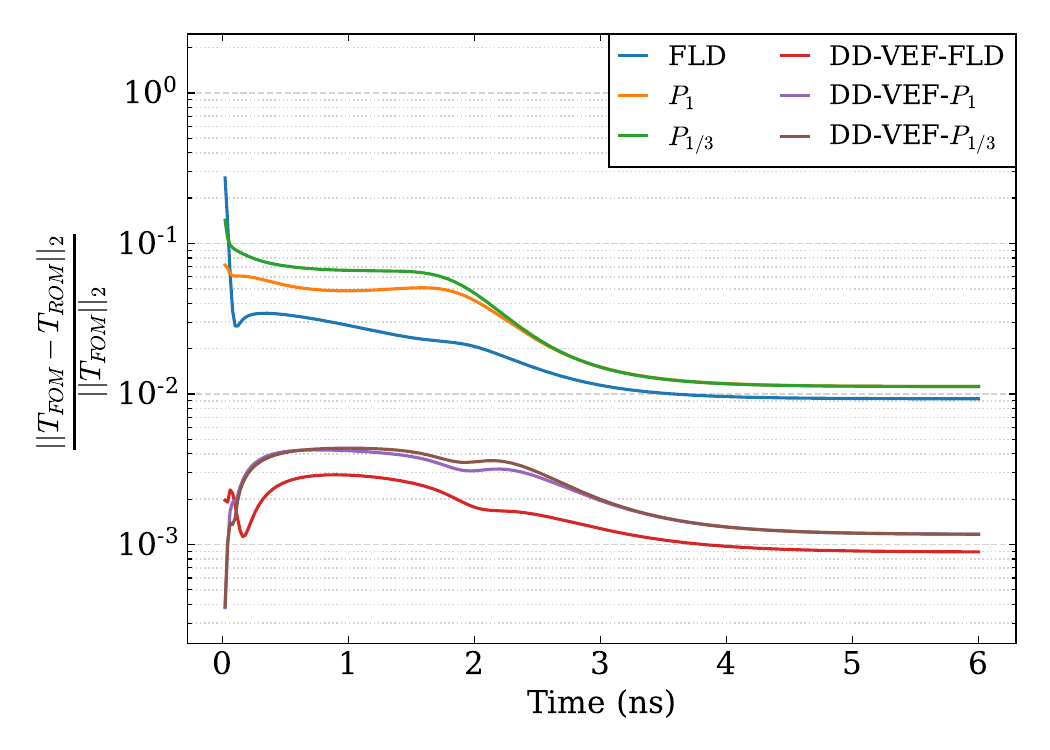}}
	\subfloat[Radiation Energy Density]{\includegraphics[width=.5\textwidth]{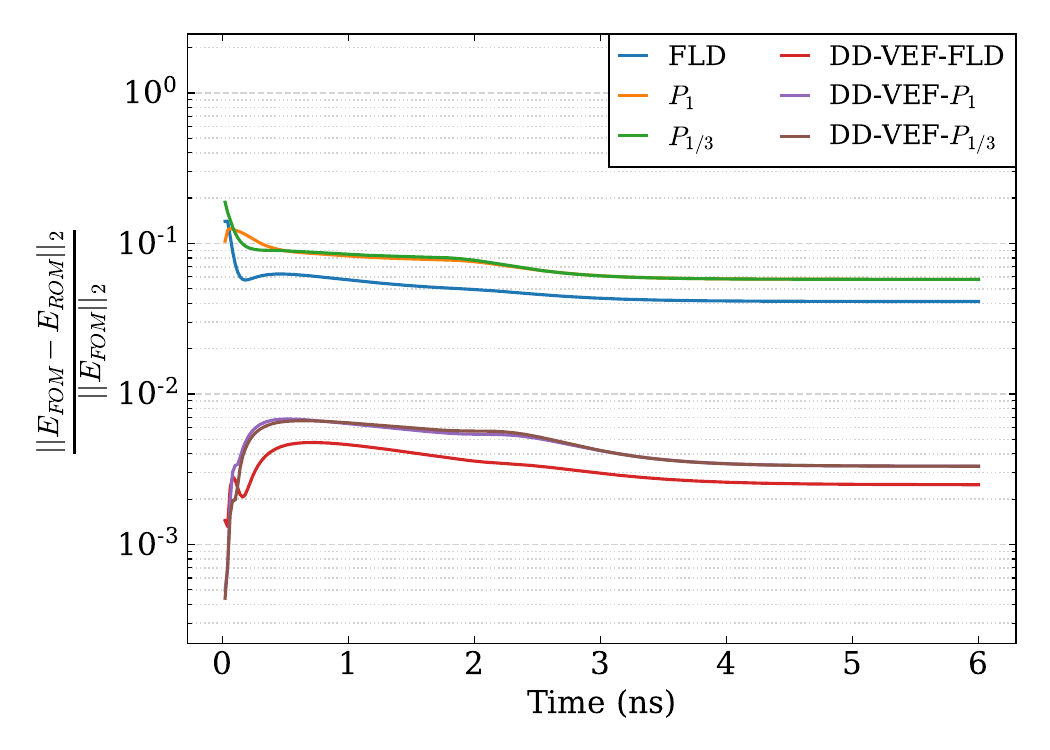}}	
	\caption{Relative errors in the 2-norm for $\tilde{T}$, $\tilde{E}$ produced by the FLD, $P_1$ and $P_{1/3}$ models, and for $T$, $E$ found with the DD-VEF model generated using each $\tilde{T}$, plotted vs time. Errors are calculated vs the FOM.
		\label{fig:TE_2nrm-errs_roms}}
\end{figure}

Figure \ref{fig:TE_2nrm-errs_roms} plots relative errors (w.r.t. the  FOM solution) for the material temperature and total radiation energy density calculated in the 2-norm over space at each instant of time in $t\in[0,6\, \text{ns}]$. Separate curves are shown for each considered diffusion model and for the DD-VEF model. In each case the DD-VEF  solution finds an increase in accuracy for $T$ and $E$ compared to the radiation diffusion solutions. The errors in $T$ and $E$ from the DD-VEF  model are on the order of $10^{-3}$ for the whole interval of time, whereas the diffusion model errors exist on order $10^{-2}$ for the majority of times. The DD-VEF  model is seen to increase the accuracy of each diffusion model by roughly an order of magnitude. The FLD model possesses the highest accuracy of all tested diffusion ROMs, and the DD-VEF  model with highest accuracy is the one applied to the FLD solution.

\begin{figure}[ht!]
	\centering
	\subfloat[$\bar{F}_R$]{\includegraphics[width=.33\textwidth]{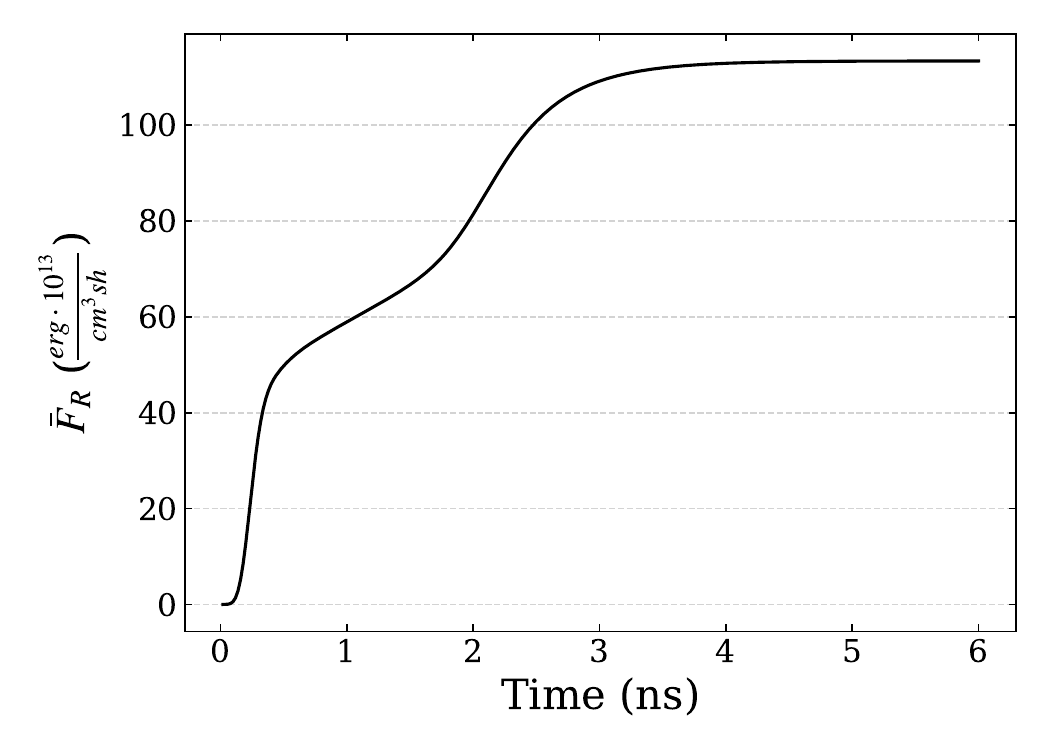}}
	\subfloat[$\bar E_R$]{\includegraphics[width=.33\textwidth]{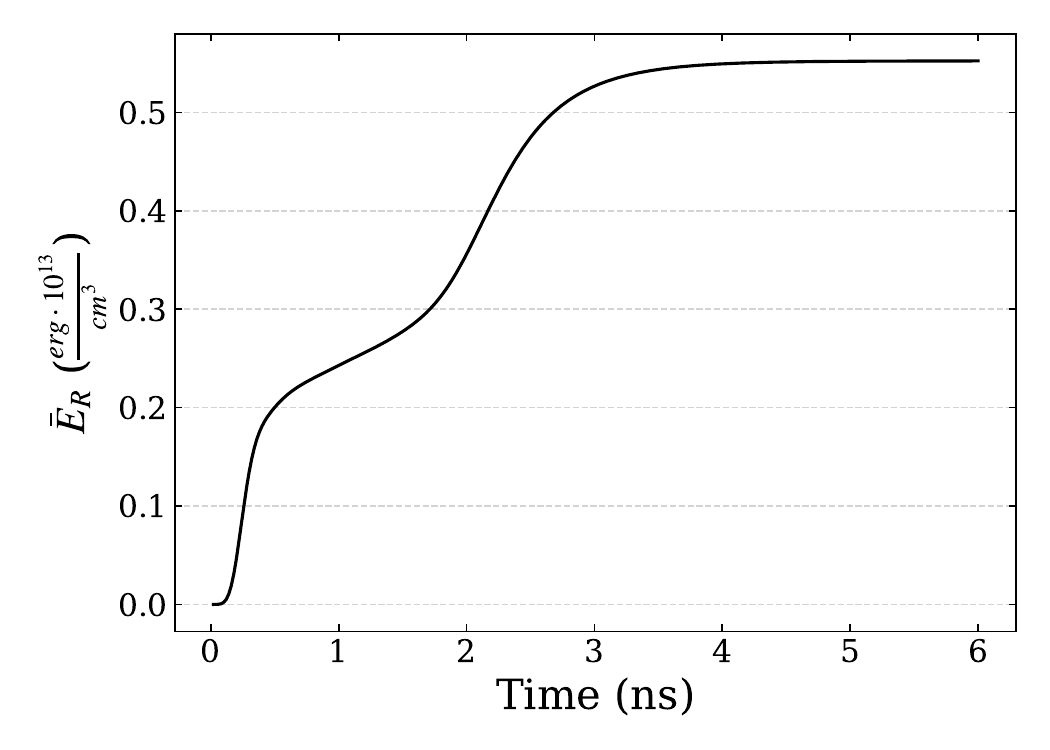}}
	\subfloat[$\bar T_R$]{\includegraphics[width=.33\textwidth]{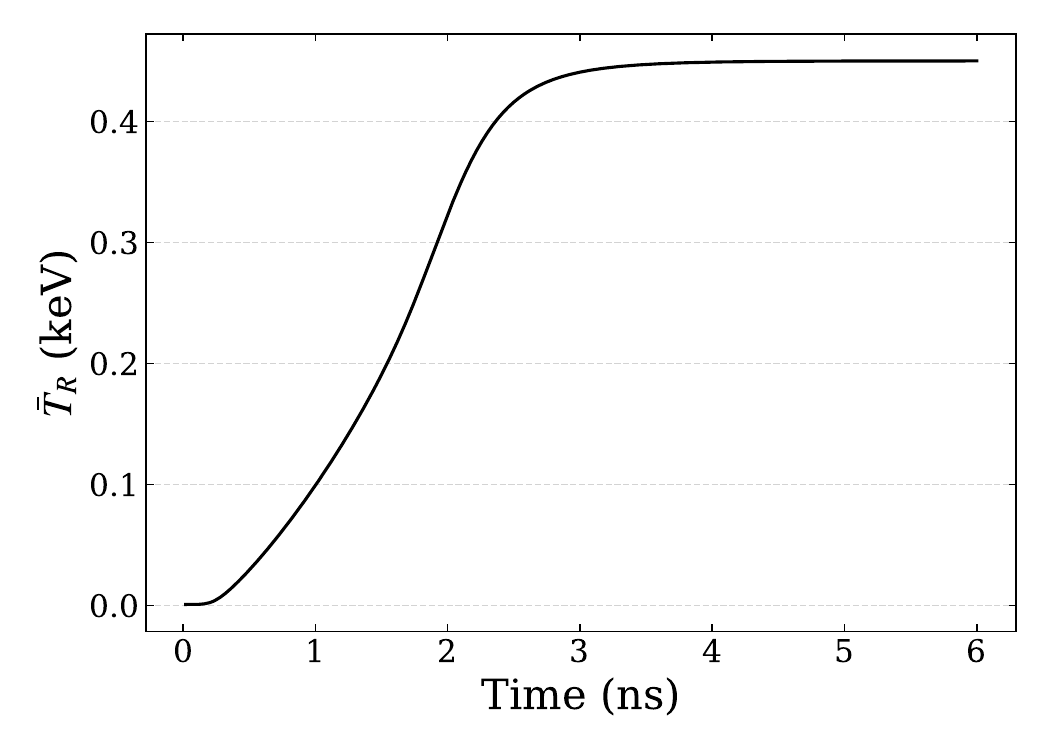}}
	\caption{FOM solution for data located at and integrated over the right boundary of the domain.
		\label{fig:rbndvals_fom} }

	\centering
	\subfloat[$\bar{F}_R$]{\includegraphics[width=.33\textwidth]{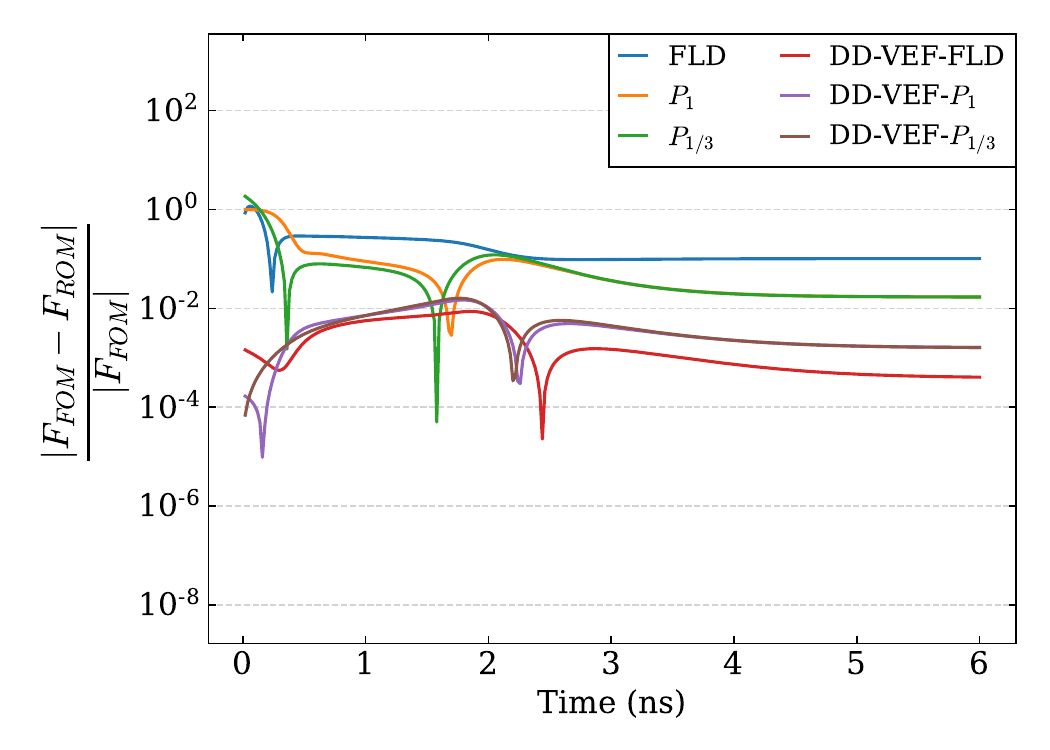}}
	\subfloat[$\bar E_R$]{\includegraphics[width=.33\textwidth]{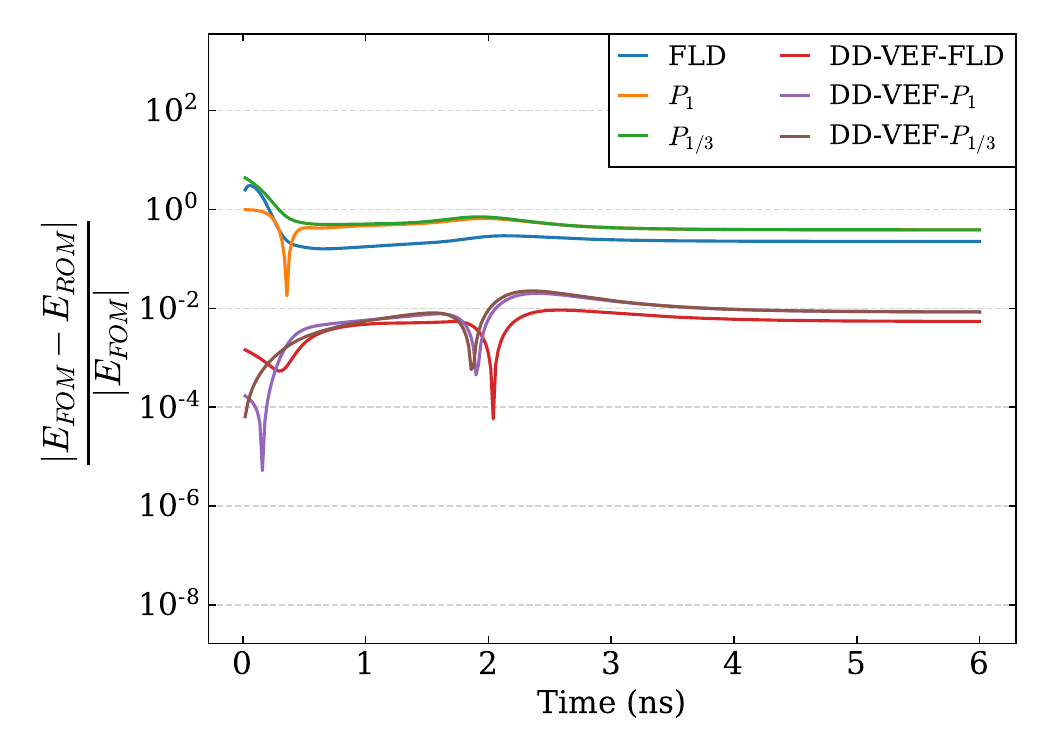}}
	\subfloat[$\bar T_R$]{\includegraphics[width=.33\textwidth]{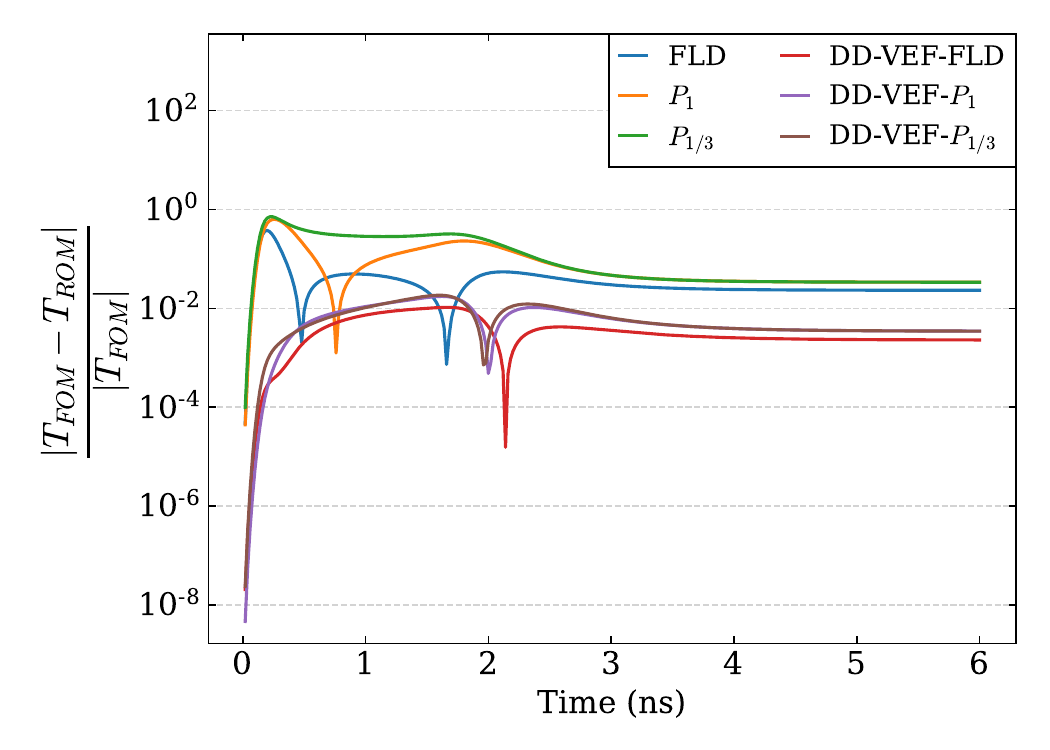}}
	\caption{Relative error for the FLD, $P_1$ and $P_{1/3}$ models, and the DD-VEF  model generated with each diffusion model solution,
		for data located at and integrated over the right boundary of the domain.
		\label{fig:rbndvals_roms} }
\end{figure}

Next we investigate the DD-VEF  model's performance in capturing the radiation wavefront as it propagates through the spatial domain.
Note that the F-C test mimics the class of supersonic radiation shock problems and experiments \cite{Moore-2015,guymer-2015,fryer-2016,fryer-2020}. One measurement of importance in these experiments concerns the time it takes for the radiation wavefront to reach the radiation-drive-opposite side of the test material \cite{Moore-2015,fryer-2016,fryer-2020}.
 A measurement of accuracy based on this wavefront-arrival metric  can be derived by comparing the TRT solution at the right boundary of the F-C test, where the radiation wavefront propagates towards. The boundary-averaged material temperatures and radiation energy density and flux are defined as
\begin{gather} \label{eq:ch4:rbnd_ints}
	\bar{F}_R = \frac{1}{L_R}\int_{0}^{L_R} \boldsymbol{e}_x\cdot\boldsymbol{F}(x_R,y)\ dy,\\[5pt]
	\bar{E}_R = \frac{1}{L_R}\int_{0}^{L_R} E(x_R,y)\ dy,\\[5pt]
	\bar{T}_R = \frac{1}{L_R}\int_{0}^{L_R} T(x_R,y)\ dy,
\end{gather}

\noindent where $L_R = x_R = 6\text{cm}$. The FOM solution for these three quantities vs time is plotted in Figure \ref{fig:rbndvals_fom}. If the DD-VEF  model can reproduce these integral quantities to acceptable levels of accuracy, they can be said to correctly reproduce the shape and propagation speed of the radiation wavefront. This is especially important to investigate given that the considered radiation diffusion models are known to produce nonphysical effects \cite{olson-auer-hall-2000,morel-2000,simmons-mihalas-2000}. Figure \ref{fig:rbndvals_roms} plots the relative error of the diffusion and DD-VEF  produced values of $\bar{F}_R$, $\bar{E}_R$ and $\bar{T}_R$. The relative errors for each quantity are decreased using the DD-VEF  model by 1-2 orders of magnitude for most instances of time. There are several time intervals where the errors for a model `spike' downwards and then come back up. These occur when there is a change of sign in the error and do not indicate that the solution is more accurate there than at other instants of time. The most dramatic increase in accuracy is for the FLD $\bar{F}_R$ by about 3 orders of magnitude. In fact, the FLD $\bar{F}_R$ is the least accurate and the DD-VEF  model $\bar{F}_R$ using the FLD solution is the most accurate of the models shown. The explanation for this effect comes from the fact that the DD-VEF  model only acts on the approximate material temperature it is given, and the FLD solution for $\bar{T}_R$ (and $T$ in general from Figure \ref{fig:TE_2nrm-errs_roms}) is the most accurate of the considered radiation diffusion models.

\begin{figure}[ht!]
	\centering
	\subfloat[$x=6$ cm, $y=3$ cm]{\includegraphics[width=.5\textwidth]{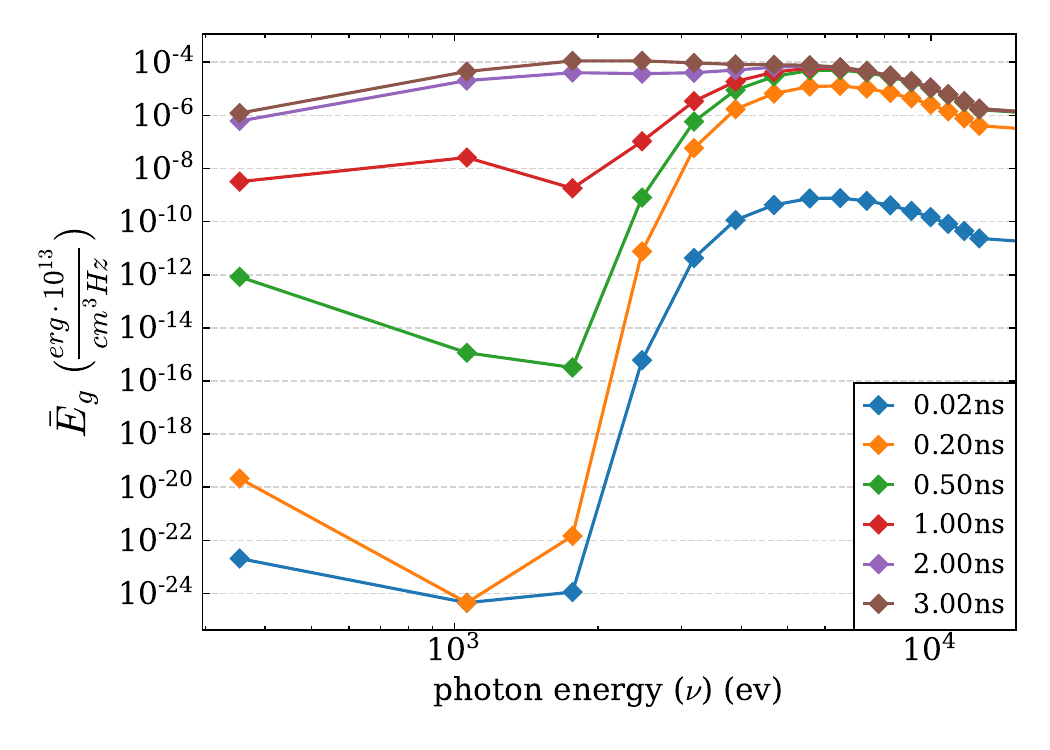}}
	\subfloat[$x=6$ cm, $y=0$ cm]{\includegraphics[width=.5\textwidth]{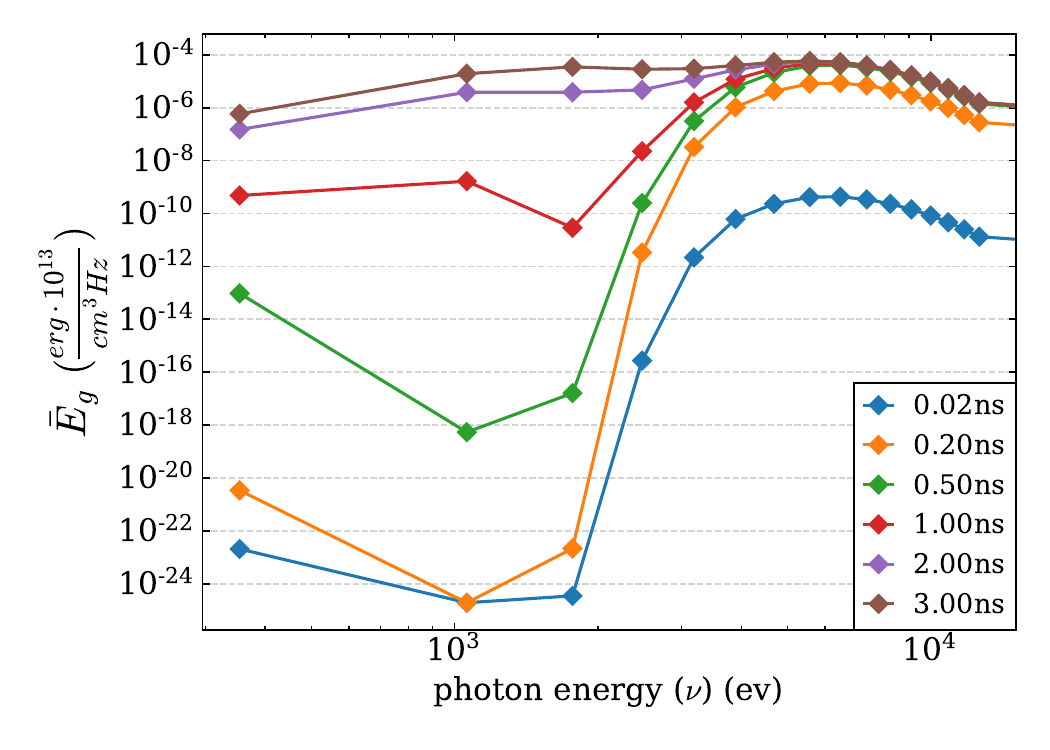}}
	\caption{Radiation energy density spectrum located at two points on the right boundary of the domain of the F-C test, taken at several time instances.
		\label{fig:Eg_spectrum} }
\end{figure}

Finally, we consider the spectrum of radiation present on the right boundary of the test domain.
Figure \ref{fig:Eg_spectrum} plots the frequency spectrum of radiation energy densities for the F-C test obtained by FOM on two points of the right boundary face ($x=6$ cm).
The spectrum of radiation present at the midpoint of the right boundary ($y=3$ cm) is shown on the left, and the radiation spectrum present at the corner of the test domain ($y=0$ cm) is displayed on the right.
Select instants of time are plotted to demonstrate how the radiation spectrum evolves.
The points on the graphs are located at the center of each discrete energy group on the frequency-axis, and the values they take on are the group-averaged radiation energy densities $\bar{E}_g={E_g}/{(\nu_{g} - \nu_{g-1})}$.
The plots have been `zoomed in' to the spectrum peak, which leaves off the final frequency group. However, this point does not deviate significantly from the position of the second to last frequency group.

\begin{figure}[ht!]
	\centering
	\subfloat[$x=6$ cm, $y=3$ cm]{\includegraphics[width=.5\textwidth]{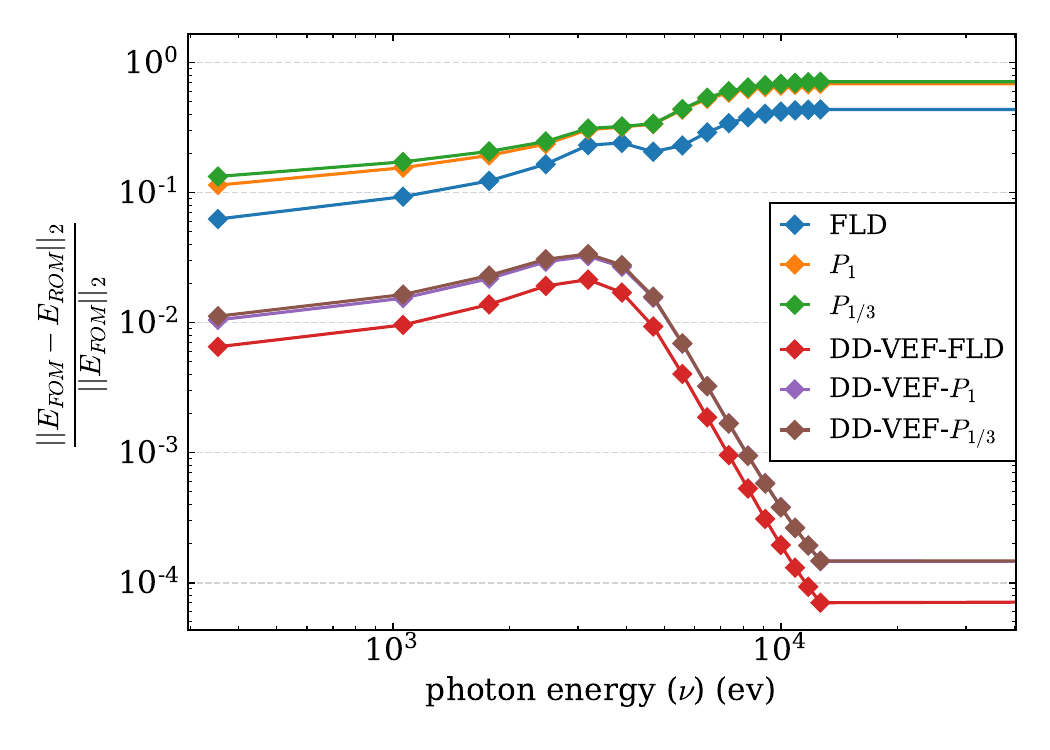}}
	\subfloat[$x=6$ cm, $y=0$ cm]{\includegraphics[width=.5\textwidth]{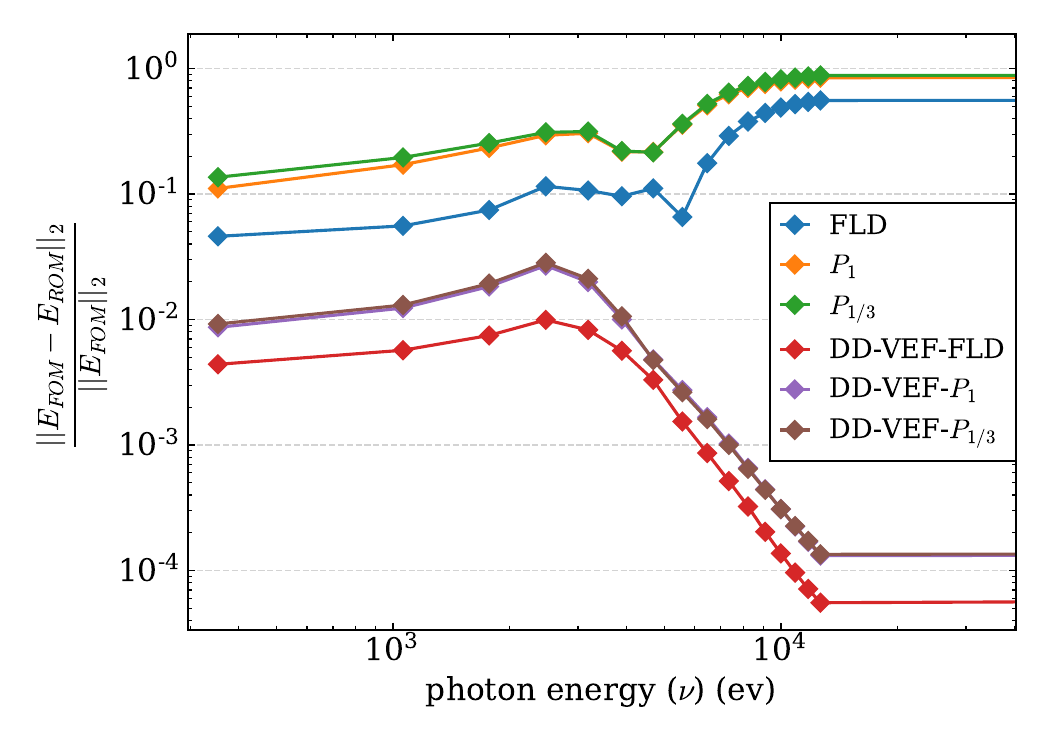}}
	\caption{Relative errors of the radiation spectrum located at the midpoint of the right boundary in the temporal 2-norm
		\label{fig:Eg_spectrum_err_2norms} }
\end{figure}

Figure \ref{fig:Eg_spectrum_err_2norms} plots the errors of each model in the radiation energy densities vs photon frequency at the midpoint of the right boundary. Errors have been collected in the relative {\it temporal} 2-norms
\begin{equation}
	\| x(t) \|_2^t=\bigg(\int_{0}^{t^\text{end}}x(t)^2dt\bigg)^{-1/2},
\end{equation}

\noindent w.r.t. the full-order solution. The DD-VEF  model is demonstrated to improve upon low-frequency group errors by roughly an order of magnitude at each considered point in space. The increase in accuracy from the diffusion solutions significantly improves as frequency increases starting from roughly $\nu=3$ KeV. This is where the peak of (non-local) radiation emanating from the boundary drive should be located in frequency, as the Planckian spectrum peaks at $\nu=2.82T$. This makes sense, since the higher frequency groups are closer to the streaming regime and should be better approximated by the transport-effects correction provided  within  the  VEF model.
The same errors calculated in the $\infty-$norm are close to those shown in the 2-norm, indicating that these results well represent the overall errors produced by these models in the radiation spectrum at all instants of time.
Note that although the last frequency group has not been included in the plots, its error value is close to that of the last shown frequency group for all models.

%
%
\section{Conclusion} \label{sec:conclusion}
In this paper a data-driven VEF  model is introduced for nonlinear TRT problems.
An approximate Eddington tensor is constructed with a transport correction method applied to radiation diffusion-based solutions to the TRT problem.
Three multigroup diffusion models were considered: a FLD model, and the $P_1$, $P_{1/3}$ models.
The DD-VEF  model provided an increase in accuracy of 1-2 orders of magnitude in the total radiation energy density and material temperature when applied to each diffusion-based solution.
The entire spectrum of radiation present at the test domain right boundary was improved upon as well.
The most significant reduction in error from the diffusion solutions in the frequency spectrum was in the high-frequency range with strong transport effects.
Possible future extensions of this DD-VEF  model include parameterization via interpolation between diffusion solutions for a series of TRT problems, or the use of other approximate models for TRT in place of  radiation-diffusion.

\section{Acknowledgements}
Los Alamos Report LA-UR-23-31255.
This research project was funded by the Sandia National Laboratory, Light Speed Grand Challenge, LDRD, Strong Shock Thrust.
This work was supported by the U.S. Department of Energy through the Los Alamos National Laboratory. Los Alamos National Laboratory is operated by Triad National Security, LLC, for the National Nuclear Security Administration of U.S. Department of Energy (Contract No. 89233218CNA000001).
The content of the information does not necessarily reflect the position or the policy of the federal government, and no official endorsement should be inferred.

\bibliographystyle{model1-num-names}
\bibliography{jmc-dya-tcdrom}

\end{document}